\documentclass[a4paper, 12pt]{amsart}

\usepackage[top=30truemm,bottom=30truemm,left=30truemm,right=30truemm]{geometry}

\usepackage{amscd, amsmath, amssymb, amsthm, graphics}

\numberwithin{equation}{section}

\theoremstyle{plain}

\newtheorem{thm}{Theorem}

\newtheorem{cor}{Corollary}
\newtheorem*{lem}{Lemma}

\theoremstyle{definition}

\newtheorem{example}{Example}
\newtheorem{rem}{Remark}

\theoremstyle{remark}

\begin{document}

\title[Weakly reflective submanifolds in isotropy irreducible spaces]
{On weakly reflective submanifolds in \\ 
compact isotropy irreducible \\
Riemannian homogeneous spaces}
\author[M. Morimoto]{Masahiro Morimoto}

\address[M. Morimoto]{Department of Mathematics, Graduate School of Science, Osaka City University. 3-3-138 Sugimoto, Sumyoshi-ku, Osaka 558-8585 Japan.}
\email{mmasahiro0408@gmail.com}

%\date{\today}



\thanks{The author was partly supported by the Grant-in-Aid for JSPS Research Fellow (No.18J14857) and by Osaka City University Advanced
Mathematical Institute (MEXT Joint Usage/Research Center on Mathematics and Theoretical Physics).}

\maketitle

\begin{abstract}
We show that for any weakly reflective submanifold of a compact isotropy irreducible Riemannian homogeneous space its inverse image under the parallel transport map is an infinite dimensional weakly reflective PF submanifold of a Hilbert space. This is an extension of the author's previous result in the case of compact irreducible Riemannian symmetric spaces. We also give a characterization of so obtained weakly reflective PF submanifolds.
\end{abstract}

\section*{Introduction}
A submanifold $N$ immersed in a Riemannian manifold $M$ is called \emph{weakly reflective} (\cite{IST09}) if for each normal vector $\xi$ at each $p \in N$ there exists an isometry $\nu_\xi$ of $M$ satisfying the conditions $\nu_\xi(p) =p$, $d \nu_\xi(\xi) = - \xi$ and $\nu_\xi(N) = N$. We call such an isometry $\nu_\xi$ a reflection with respect to $\xi$. If every $\nu_\xi$ can be chosen from a particular subgroup $S$ of the isometry group $I(M)$ then we call $N$ \emph{$S$-weakly reflective}. By definition weakly reflective submanifolds are \emph{austere} (\cite{HL82}): for each normal vector $\xi$ the set of eigenvalues with multiplicities of the shape operator $A_\xi$ is invariant under the multiplication by $(-1)$. Thus weakly reflective submanifolds are minimal submanifolds.  A singular orbit of a cohomogeneity one action is a typical example of a weakly reflective submanifold (\cite{Pod97}, \cite{IST09}). It is an interesting problem to determine weakly reflective orbits in an isometric action of a Lie group (e.g.\ \cite{IST09}, \cite{Eno18}).

Recently the author \cite{M1}  introduced the concept of weakly reflective submanifolds into a class of proper Fredholm (PF) submanifolds in Hilbert spaces (Terng \cite{Ter89}) and 
showed that if $N$ is a weakly reflective submanifold of an irreducible Riemannian symmetric space $G/K$ of compact type then its inverse image $\Phi_K^{-1}(N)$ under the parallel transport map $\Phi_K: V_\mathfrak{g} \rightarrow G/K$ (\cite{TT95}) is an infinite dimensional weakly reflective PF submanifold of the Hilbert space $V_\mathfrak{g} := L^2([0,1], \mathfrak{g})$ consisting of all $L^2$-paths from $[0,1]$ to the Lie algebra $\mathfrak{g}$ of $G$ (\cite[Theorem 8]{M1}). Using this result many examples of infinite dimensional weakly reflective PF submanifolds were obtained from finite dimensional weakly reflective submanifolds in $G/K$. The purpose of this paper is to extend that result to the case $G/K$ is a compact isotropy irreducible Riemannian homogeneous space and to give a characterization of so obtained weakly reflective PF submanifolds in terms of $S$-weak reflectivity defined above. The results are summarized in Section 2. Their proofs are given in Section 3. To do these, in Section 1 we prepare the setting of $P(G,H)$-actions and the parallel transport map over a compact Lie group $G$ which is \emph{not} necessarily connected. This generalized setting is important in the formulation of our results. 

\section{Preliminaries}

Let $G$ be a compact Lie group with Lie algebra $\mathfrak{g}$. Denote by $G_0$ its identity component. Choose an $\operatorname{Ad}(G)$-invariant inner product of $\mathfrak{g}$ and equip the corresponding bi-invariant inner product with $G$. Denote by $\mathcal{G} := H^1([0,1], G)$ the Hilbert Lie group of all Sobolev $H^1$-paths from $[0,1]$ to $G$ and by $V_\mathfrak{g} := L^2([0,1], \mathfrak{g})$ the Hilbert space of all $L^2$-paths from $[0,1]$ to $\mathfrak{g}$. For each $a \in G$ (resp.\ $x \in \mathfrak{g}$) we denote by $\hat{a} \in \mathcal{G}$ (resp.\ $\hat{x} \in V_\mathfrak{g}$) the constant path which values at $a$ (resp.\ $x$). Define the $\mathcal{G}$-action on $V_\mathfrak{g}$ via the gauge transformations:
\begin{equation}\label{gauge}
g * u := gug^{-1} - g' g^{-1}
, \quad 
g \in \mathcal{G}, \ u \in V_\mathfrak{g},
\end{equation}
where $g'$ denotes the weak derivative of $g$. We know that this action is isometric, transitive, proper and Fredholm    (\cite[Theorem 5.8.1]{PT88}).

Let $H$ be a closed subgroup of $G \times G$. Then $H$ acts on $G$ isometrically by
\begin{equation} \label{Haction}
(b_1, b_2) \cdot a := b_1 a b_2^{-1}
,\quad
a, b_1, b_2 \in G.
\end{equation}
Define a Lie subgroup $P(G,H)$ of $\mathcal{G}$ by 
\begin{equation*}
P(G, H) := \{g \in \mathcal{G} \mid (g(0), g(1)) \in H\}.
\end{equation*}
The induced action of $P(G,H)$ on $V_\mathfrak{g}$ is called the \emph{$P(G,H)$-action}. We know that the $P(G,H)$-action is isometric, proper and Fredholm (\cite[p.\ 132]{Ter95}). 

The \emph{parallel transport map} (\cite{KT93}, \cite{Ter95}) $\Phi: V_\mathfrak{g} \rightarrow G_0$ is a Riemannian submersion defined by 
\begin{equation*}
\Phi(u) :=E_u(1), \quad u \in V_\mathfrak{g},
\end{equation*}
where $E_u \in \mathcal{G}$ is the unique solution to the linear ordinary differential equation
\begin{equation*}
\left\{
\begin{array}{l}
E_u^{-1} E_u' = u, 
\\
E_u(0) = e.
\end{array}
\right.
\end{equation*}
Let $H$ be a closed subgroup of $G_0 \times G_0$. Then it follows (\cite[Proposition 1.1]{Ter95}) that for each $g \in P(G_0,H)$ and $u \in V_\mathfrak{g}$
\begin{equation} \label{equiv4}
\Phi(g * u) = (g(0), g(1)) \cdot \Phi(u)
 \quad \text{and} \quad
P(G_0, H) * u = \Phi^{-1}(H \cdot \Phi(u)).
\end{equation}
In general, if $N$ is a closed submanifold of $G_0$ then the inverse image $\Phi^{-1}(N)$ is a PF submanifold of $V_\mathfrak{g}$ (\cite[Lemma 5.1]{TT95}).

Let $K$ be a closed subgroup of $G_0$ with Lie algebra $\mathfrak{k}$. Denote by $\mathfrak{g}= \mathfrak{k} + \mathfrak{m}$ the orthogonal direct sum  decomposition. Restricting the $\operatorname{Ad}(G)$-invariant inner product of $\mathfrak{g}$ to $\mathfrak{m}$ we equip the induced $G_0$-invariant Riemannian metric with the homogeneous space $G_0/K$. We call such a metric the normal homogeneous metric and $G_0/K$ the  compact normal homogeneous space. Denote by $\pi : G_0 \rightarrow G_0/K$ the natural projection, which is a Riemannian submersion with totally geodesic fiber. The \emph{parallel transport map $\Phi_K$ over $G_0/K$} is a Riemannian submersion defined by 
\begin{equation} \label{ptm2} 
\Phi_K := \pi \circ \Phi : V_\mathfrak{g} \rightarrow G_0 \rightarrow G_0/K.
\end{equation}

\section{Results}

Let $M$ be a compact Riemannian homogeneous space, that is, the group of isometries $G := I(M)$ is compact and the identity component $G_0 := I_0(M)$ acts transitively on $M$. Fix $p \in M$ and denote by $K := I_0(M)_p$ the isotropy subgroup of $G_0$ at $p \in M$. Then we have a diffeomorphism $G_0/K \cong M$. Choose an inner product of the Lie algebra $\mathfrak{g}$ which is invariant not only under $\operatorname{Ad}(G_0)$ but also under $\operatorname{Ad}(G)$. We equip the corresponding bi-invariant metric with $G_0$ and the normal homogeneous metric with $G_0/K$. Here note that $M$ and $G_0/K$ are \emph{not} isometric in general. To avoid this we suppose that $M$ is \emph{isotropy irreducible}: the linear isotropy representation of $K$ on the tangent space $T_p M$ is $\mathbb{R}$-irreducible. From this condition the $G_0$-invariant Riemannian metric on $G_0/K$ is unique up to scaling and thus we can (and will) always choose an $\operatorname{Ad}(G)$-invariant inner product of $\mathfrak{g}$ so that the normal homogeneous space $G_0/K$ is isometric to $M$. For the geometry, structure and the classification of isotropy irreducible  homogeneous spaces, see \cite{WZ91}, \cite{Wol68} and \cite{Wol84}.

The following theorem is the main result of this paper:
\begin{thm}\label{thm1}
Let $M$ be a compact isotropy irreducible Riemannian homogeneous space and $N$ a closed submanifold of $M$. Denote by $G$ the group of isometries on $M$, by $G_0$ its identity component, by $K$ the isotropy subgroup of  $G_0$ at a fixed $p \in M$ and by $\Phi_K = \pi \circ \Phi: V_{\mathfrak{g}} \rightarrow G_0 \rightarrow G_0/K = M$ the parallel transport map. Then for a closed  subgroup $S$ of $G$ satisfying $a S a^{-1} = S$ for all $a \in G_0$ the following are equivalent: 
\begin{enumerate}
\item $N$ is an $S$-weakly reflective submanifold of $M$,
\item $\pi^{-1}(N)$ is an $(S \times S_p)$-weakly reflective submanifold of $G_0$,  
\item $\Phi_K^{-1}(N)$ is a $P(S, S \times S_{p})$-weakly reflective PF submanifold of $V_\mathfrak{g}$,
\end{enumerate}
where $S_p$ denotes the isotropy subgroup of $S$ at $p$.
\end{thm}

\begin{rem} Strictly speaking ``$(S \times S_p)$-weakly reflective'' in the statement (ii) above should be written by ``$((S \times S_p) \cap I(G_0))$-weakly reflective'' because not all  elements of $(S \times S_p)$ preserve $G_0$ and induce isometries on $G_0$. However we will continue to use such an abbreviation for simplicity. 
\end{rem}

\begin{rem} In the previous paper \cite[Theorem 8]{M1} $M$ was assumed to be a symmetric space of compact type and only the case $S = G$ was considered. Theorem \ref{thm1} here does not require such assumptions and moreover it characterizes the weakly reflective PF submanifold $\Phi_K^{-1}(N)$ precisely. \end{rem}

Considering the case $S = G$ we obtain:
\begin{cor}\label{cor1}
Let $M = G_0/K$ and $N$ be as above. Then the following are equivalent:
\begin{enumerate}
\item $N$ is a weakly reflective submanifold of $M$,
\item $\pi^{-1}(N)$ is a $(G \times G_p)$-weakly reflective submanifold of $G_0$,  
\item $\Phi_K^{-1}(N)$ is a $P(G, G \times G_p)$-weakly reflective PF submanifold of $V_{\mathfrak{g}}$. 
\end{enumerate}
\end{cor}

Considering the case $S = G_0$ we obtain:
\begin{cor}\label{cor2}
Let $M = G_0/K$ and $N$ be as above. Then the  following are equivalent:
\begin{enumerate}
\item $N$ is a $G_0$-weakly reflective submanifold of $M$,
\item $\pi^{-1}(N)$ is a $(G_0 \times K)$-weakly reflective submanifold of $G_0$,  
\item $\Phi_K^{-1}(N)$ is a $P(G_0, G_0 \times K)$-weakly reflective PF submanifold of $V_\mathfrak{g}$. 
\end{enumerate}
\end{cor}

Next we suppose a certain homogeneous condition for $N$ and the condition $S \subset G_0$. However we do not suppose that $S$ satisfies $aS a^{-1} = S$ for \emph{all} $a \in G_0$:
\begin{thm} \label{thm2}
Let $M = G_0 /K$ and $N$ be as in Theorem $\ref{thm1}$. Suppose that $N$ is an orbit of a closed subgroup $U$ of $G_0$. Then for a closed subgroup $S$ of $G_0$ satisfying $aSa^{-1} = S$ for all $a \in U$ the following are equivalent:
\begin{enumerate}
\item $N$ is an $S$-weakly reflective submanifold of $M$,
\item $\pi^{-1}(N)$ is an $(S \times K)$-weakly reflective submanifold of $G_0$,  
\item $\Phi_K^{-1}(N)$ is a $P(G_0, S \times K)$-weakly reflective PF submanifold of $V_\mathfrak{g}$.
\end{enumerate}
\end{thm}

Considering the case $S = U$ we obtain:
\begin{cor}\label{cor3}
Let $M = G_0/K$, $N$ and $U$ be as in Theorem \ref{thm2}. Then the following are equivalent:
\begin{enumerate}
\item $N$ is an $U$-weakly reflective submanifold of $M$,
\item $\pi^{-1}(N)$ is an $(U \times K)$-weakly reflective submanifold of $G_0$,  
\item $\Phi_K^{-1}(N)$ is a $P(G_0, U \times K)$-weakly reflective PF submanifold of $V_\mathfrak{g}$.
\end{enumerate}
\end{cor}

Considering the case $S = U_0$ we obtain:
\begin{cor}\label{cor4}
Let $M = G_0/K$, $N$ and $U$ be as in Theorem \ref{thm2}. Then the following are equivalent:
\begin{enumerate}
\item $N$ is an $U_0$-weakly reflective submanifold of $M$,
\item $\pi^{-1}(N)$ is an $(U_0 \times K)$-weakly reflective submanifold of $G_0$,  
\item $\Phi_K^{-1}(N)$ is a $P(G_0, U_0 \times K)$-weakly reflective PF submanifold of $V_\mathfrak{g}$.
\end{enumerate}
\end{cor}

Note that if $N$ is an orbit $U \cdot aK$ through $aK \in G_0/K$ then we have $\pi^{-1}(N) = (U \times K) \cdot a$ and $\Phi_K^{-1}(N) = P(G_0, U \times K) * u$ for $u \in \Phi^{-1}(a)$ by \eqref{equiv4}. Applying above results to examples of weakly reflective submanifolds in $M$ we obtain examples of infinite dimensional weakly reflective PF submanifolds in Hilbert spaces as follows:
\begin{example}
Let $M = G_0/K$ be as above and $U$ a closed subgroup of $G_0$. Suppose that the $U$-action on $M$ is of cohomogeneity one. Then any singular orbit $N = U \cdot aK$ is a $U$-weakly reflective submanifold of $M$ (\cite{Pod97}, \cite{IST09}). Applying Corollary  \ref{cor3} to $N$ the orbit $P(G_0, U \times K) * u$ through $u \in \Phi^{-1}(a)$ is a $P(G_0, U \times K)$-weakly reflective PF submanifold of $V_\mathfrak{g}$. In particular if $\operatorname{codim} N \geq 2$ it was essentially shown that $N$ is a $U_0$-weakly reflective submanifold of $M$. Then by Corollary \ref{cor4} the orbit  $P(G_0, U \times K) * u$ is a $P(G_0, U_0 \times K)$-weakly reflective PF submanifold of $V_\mathfrak{g}$. 
\end{example}

\begin{example}
Let $n \in \mathbb{Z}_{\geq 2}$. Set  $M := S^{2n-1}(\sqrt{2}) \subset \mathbb{R}^{2n}$ and $N := S^{n-1}(1) \times S^{n-1}(1) \subset M$. It was shown (\cite{IST09}) that $N$ is an $O(2n)$-weakly reflective submanifold of $M$. Fix $p \in N$. Set $(G_0,K) := (SO(2n), SO(2n)_p)$ and $U := SO(n) \times SO(n)$ so that $M = G_0/K$ and $N = U \cdot eK$. Applying Corollary \ref{cor1} to $N$ the orbit $P(G_0, U \times K) * \hat{0}$ is a $P(O(2n), O(2n) \times O(2n)_p)$-weakly reflective PF submanifold of $V_\mathfrak{g}$. In particular if $n$ is even then it was essentially shown that $N$ is an $SO(2n)$-weakly reflective submanifold of $M$. Then by Corollary \ref{cor2} the orbit $P(G_0, U \times K) * \hat{0}$ is a $P(SO(2n), SO(2n) \times SO(2n)_p)$-weakly reflective PF submanifold of $V_\mathfrak{g}$.
\end{example}

\begin{example}
Consider the symmetric pair $(SO(7), SO(3) \times SO(4))$ and the cohomogeneity one action of the exceptional Lie group $G_2$ on the symmetric space $M = SO(7)/(SO(3) \times SO(4))$. Then there exists a unique weakly reflective principal orbit $N = G_2 \cdot aK$, which was shown to be $O(7)$-weakly reflective (\cite{Eno18}). Applying Theorem \ref{thm1} to $N$ the orbit $P(SO(7), G_2 \times SO(3) \times SO(4)) * u$ through $u \in \Phi^{-1}(a)$ is a $P(O(7), O(7) \times O(3) \times O(4))$-weakly reflective PF submanifold of $V_\mathfrak{g}$.
\end{example}

\begin{rem}
Recently Taketomi \cite{Tak18} introduced a generalized concept of weakly reflective submanifolds, namely arid submanifolds. The theorems and corollaries in this section are still valid in the arid case (see also \cite{M2}). 
\end{rem}

\section{Proofs of the theorems}
To prove the theorems we first recall several facts. Let $G$ be a compact Lie group with a bi-invariant Riemannian metric and $\Phi: V_\mathfrak{g} \rightarrow G_0$ the parallel transport map. We know that the $P(G_0, G_0 \times \{e\})$-action on $V_\mathfrak{g}$ is simply transitive (\cite[Corollary 4.2]{TT95}). By \eqref{equiv4} the following diagram commutes for each $g \in P(G_0, G_0 \times \{e\})$:
\begin{equation} \label{commute1}
\begin{CD}
V_\mathfrak{g}@>g* >> V_\mathfrak{g}
\\
@V\Phi VV @V\Phi VV
\\
G_0 @> (g(0),\, e) >>  G_0.\!
\end{CD}
\end{equation}
Let $K$ be a closed subgroup of $G_0$ and $G_0/K$ the compact normal homogeneous space with projection $\pi: G_0 \rightarrow G_0/K$. For each $a \in G_0$ we denote by $l_a$ the left translation by $a$ and $L_a$ the  isometry on $G_0/K$ defined by $L_a(bK) := (ab)K$. Then we have the commutative diagram:
\begin{equation}\label{commute2}
\begin{CD}
G_0 @>l_a>> G_0
\\
@V\pi VV @V\pi VV
\\
G_0/K @>L_a>> G_0/K.\!
\end{CD}
\end{equation}
The following lemma was shown in \cite[Lemma 2]{M1}. 
\begin{lem}[\cite{M1}] \label{weaksubm}
Let $\mathcal{M}$ and $\mathcal{B}$ be Riemannian Hilbert manifolds, $\phi: \mathcal{M} \rightarrow \mathcal{B}$ a Riemannian submersion and $N$ a closed submanifold of $\mathcal{B}$. Fix $p \in \phi^{-1}(N)$ and $X \in T^\perp_p \phi^{-1}(N)$. Suppose that $\nu_{\mathcal{M}}$ is an isometry of $\mathcal{M}$ fixing $p$, that $\nu_{\mathcal{B}}$ is an isometry of $\mathcal{B}$ fixing $\phi(p)$ and that the diagram
\begin{equation*}
\begin{CD}
\mathcal{M} @>\nu_{\mathcal{M}}>> \mathcal{M}
\\
@V\phi VV @V\phi VV
\\
\mathcal{B} @>\nu_{\mathcal{B}}>> \mathcal{B}
\end{CD}
\end{equation*}
commutes. Then the following are equivalent:
\begin{enumerate}
\item[\textup{(i)}] $\nu_\mathcal{M}$ satisfies $\nu_\mathcal{M}(\phi^{-1}(N)) = \phi^{-1}(N)$ and $d\nu_{\mathcal{M}}(X) = - X$,
\item[\textup{(ii)}] $\nu_\mathcal{B}$ satisfies $\nu_\mathcal{B}(N) = N$ and $d \nu_{\mathcal{B}}(d \phi(X)) = - d \phi (X)$. 
\end{enumerate} 
\end{lem}

We are now in position to prove the theorems.

\medskip

\noindent\textbf{Proof of Theorem \ref{thm1}.}
By identification $M = G_0/K$ we have  $p=eK$. From the assumption, for each $a \in G_0$ and each $g \in P(G_0, G_0 \times \{e\})$ we have 
\begin{align*} 
&a^{-1}S_{aK} a= S_{eK}, \quad
(a,e)^{-1}(S \times S_{eK})_{a}(a,e) 
= 
(S \times S_{eK})_{e} = \Delta(S_{eK})
\quad \text{and} 
\\ 
&g^{-1}P(S, S \times S_{eK})_{g* \hat{0}}g 
=
P(S, S \times S_{eK})_{\hat{0}}
=
\{ \hat{s} \in \hat{S} \mid s \in S_{eK}\}.
\end{align*}
 ``(i) $\Rightarrow$ (ii)'':  Let $a \in \pi^{-1}(N)$ and $w \in T_a^\perp \pi^{-1}(N) $. Set $\eta:= d \pi (w) \in T^\perp_{aK} N $, $N' := L_a^{-1}(N)$ and $\xi := d L_a^{-1}(\eta) \in T^\perp_{eK}N'$. Denote by $v \in T_e^\perp \pi^{-1}(N')$ the horizontal lift of $\xi$. From (\ref{commute2}) we have $l_a (\pi^{-1}(N')) = \pi^{-1}(N)$ and $d l_a (v) = w$. Thus to show the existence of a reflection $\nu_w$ with respect to $w$ satisfying $\nu_w \in (S \times S_{eK})_a$ it suffices to construct an reflection $\nu_v$ with respect to $v$ satisfying $\nu_v \in (S \times S_{eK})_e$. Take a reflection $\nu_\eta$ with respect to $\eta$ satisfying $\nu_\eta \in S_{aK}$. Define a reflection $\nu_{\xi}$ with respect to $\xi$ by $\nu_\xi := L_a^{-1} \circ \nu_\eta \circ L_a$. Then $\nu_\xi \in S_{eK}$. Since $G_0$ is a normal subgroup of $G$ an automorphism $\nu_v: G_0 \rightarrow G_0$, $b \mapsto \nu_v(b)$ is well-defined by
\begin{equation*}
\nu_{v}(b)
:= 
\nu_\xi \circ b  \circ \nu_\xi^{-1},
\end{equation*}
that is,
\begin{equation*}
\nu_v = \operatorname{Ad}^G(\nu_\xi)|_{G_0}.
\end{equation*}
Clearly $\nu_v \in (S \times S_{eK})_e$. Note that $\nu_v$ is an isometry of $G_0$ because 
we fixed a bi-invariant Riemannian metric on $G_0$ coming from an $\operatorname{Ad}(G)$-invariant inner product of $\mathfrak{g}$. Moreover we have
\begin{align*}
(\nu_\xi \circ \pi) (b) 
&=
\nu_\xi(bK)
=
(\nu_\xi \circ b) (eK)
=
(\nu_\xi \circ b  \circ \nu_\xi^{-1})(eK) \quad \text{and}
\\
(\pi \circ \nu_v)(b)
&=
\nu_v(b)K
=
\nu_v(b) (eK)
=
(\nu_\xi \circ b \circ \nu_\xi^{-1})(eK)
\end{align*}
for all $b \in G_0$.  This shows that the diagram
\begin{equation*}
\begin{CD}
G_0 @>\nu_v>> G_0
\\
@V \pi VV @V \pi VV
\\
\ G_0/K\  @>\nu_{\xi}>> \ G_0/K\ 
\end{CD}
\end{equation*}
commutes. Thus by Lemma,  $\nu_v$ is a reflection with respect to $v$. This shows (ii).
 
``(ii) $\Rightarrow$ (i)'': Let $aK \in N$ and $\eta \in T_{aK}^\perp N$. Denote by $w \in T^\perp_a \pi^{-1}(N)$ the horizontal lift of $\eta$. Define $N'$, $\xi$, $v$ as above. Take a reflection $\nu_{w}$ with respect to $w$  satisfying $\nu_w \in (S \times S_{eK})_a$. Then a reflection with respect to $v$ is defined by $\nu_v := (a,e)^{-1} \circ \nu_w \circ (a,e)$ so that $\nu_v \in (S \times S_{eK})_e$. Thus there exists $s \in S_{eK}$ such that $\nu_v = (s,s)$. Define an isometry $\nu_\xi$ of $G_0/K$ by $\nu_\xi := s$. Then by Lemma, $\nu_\xi$ is a reflection with respect to $\xi$ satisfying $\nu_\xi \in S_{eK}$. Hence a reflection $\nu_\eta$ with respect to $\eta$ is defined by $\nu_\eta := l_a \circ \nu_\xi \circ l_a^{-1}$ so that $\nu_\eta \in S_{aK}$ and (i) follows.

``(ii) $\Rightarrow$ (iii)'': Set $Q := \pi^{-1}(N)$ so that $\Phi_K^{-1}(N) = \Phi^{-1}(Q)$. Let $u \in \Phi^{-1}(Q)$ and $X \in T_u^\perp \Phi^{-1}(Q) $. Take $g \in P(G_0, G_0 \times \{e\})$ so that $u = g * \hat{0}$. Set $a := \Phi(u) = g(0)$, $w := d \Phi(X) \in T_{a} ^\perp Q $, $Q' := a^{-1}Q$ and $v := dl_a^{-1}( w) \in T_e^\perp Q'$. Denote by $\hat{v} \in T^\perp_{\hat{0}} \Phi^{-1}(Q')$ the horizontal lift of $v$. From (\ref{commute1}) we have $g * (\Phi^{-1}(Q')) = \Phi^{-1}(Q)$ and $(d g*) \hat{v} = X$. Thus to show the existence of a reflection $\nu_X$ with respect to $X$ satisfying $\nu_X \in P(S, S \times S_{eK})_u$ it suffices to construct a reflection $\nu_{\hat{v}}$ with respect to $\hat{v}$ satisfying $\nu_{\hat{v}} \in P(S, S \times S_{eK})_{\hat{0}}$. Take a reflection $\nu_{w}$ with respect to $w$ satisfying $\nu_w \in (S \times S_{eK})_a$. Then a reflection  $\nu_{v}$ with respect to $v$ is defined  by 
$
\nu_{v} := (a, e)^{-1} \circ \nu_{w} \circ (a,e)
$
so that $\nu_v \in (S \times S_{eK})_e$. Thus there exists $s \in S_{eK}$ so that $\nu_v = (s,s)$. We define a linear orthogonal transformation $\nu_{\hat{v}}$ of $V_\mathfrak{g}$  by
\begin{equation*}
\nu_{\hat{v}} (u) 
:= 
d \nu_v \circ u
=
s u s^{-1}
=
\hat{s} * u,
\quad 
u \in {V_\mathfrak{g}}.
\end{equation*}
Clearly $\nu_{\hat{v}} \in P(S, S \times S_{eK})_{\hat{0}}$. Moreover since $\nu_v$ is an automorphism of $G_0$ we have $\nu_{\hat{v}}(g* \hat{0}) = (\nu_{v} \circ g) * \hat{0}$ for all $g \in P(G_0, G_0 \times G_0)$. This together with (\ref{equiv4}) implies that the following diagram commutes:
\begin{equation*}
\begin{CD}
V_\mathfrak{g} @>\nu_{\hat{v}}>> V_\mathfrak{g}
\\
@V \Phi VV @V \Phi VV
\\
\ G_0\  @>\nu_v >> \ G_0.
\end{CD}
\end{equation*}
From Lemma  $\nu_{\hat{v}}$ is a refection with respect to $\hat{v}$. This shows (iii).

``(iii) $\Rightarrow$ (ii)'': Let $a \in Q$, $w \in T_a^\perp Q$ and $u \in \Phi^{-1}(a)$. Denote by $X \in T^\perp_u \Phi^{-1}(Q)$ the horizontal lift of $w$. Take $g \in P(G_0, G_0 \times \{e\})$ so that $g*\hat{0} = u$ and define  $Q'$, $v$, $\hat{v}$ as above.  Take a reflection $\nu_X$ with respect to $X$ satisfying $\nu_X \in P(S, S \times S_{eK})_u$. Then a reflection $\nu_{\hat{v}}$ with respect to $\hat{v} \in T^\perp_{\hat{0}} \Phi^{-1}(Q')$ is defined by $\nu_{\hat{v}} := (g*)^{-1} \circ \nu_X \circ (g*)$ so that $\nu_v \in P(S, S \times S_{eK})_{\hat{0}}$. Thus there exists $s \in S_{eK}$ such that $\nu_{\hat{v}}= \hat{s}$. Define an isometry $\nu_v$ of $G_0$ by $\nu_v :=(s,s)$. From Lemma, $\nu_v$ is a reflection with respect to $v$ satisfying $\nu_v \in (S \times S_{eK})_e$. Thus a reflection $\nu_w$  with respect to $w$ is defined by $\nu_{w} := l_a \circ \nu_v \circ l_a^{-1}$ so that $\nu_w \in (S \times S_{eK})_a$. This shows (ii). \qed

\medskip

\noindent\textbf{Proof of Theorem \ref{thm2}.}
Take $b \in \pi^{-1}(N)$ and $u \in \Phi^{-1}(b)$. Then $N = U \cdot bK$, $\pi^{-1}(N) = (U \times K) \cdot b$ and $\Phi_K^{-1}(N) = P(G_0, U \times K)*u$. Choose $g \in P(G_0 , G_0 \times \{e\})$ so that $u = g * \hat{0}$. Then $b = \Phi(u) = g(0)$. Set $U' := b^{-1} U b$. By \eqref{commute1} we have
\begin{align*}
&
L_b(U' \cdot eK) = U \cdot bK,
\quad
l_b((U'\times K) \cdot e)
=
(U \times K) \cdot b
\quad \text{and}
\\
& 
g* (P(G_0, U' \times K)*\hat{0})
=
P(G_0, U \times K) * u. 
\end{align*}
Set $S' := b^{-1}Sb$. Then $a' S' (a')^{-1} = S'$ for all $a' \in U'$. By definition we have
\begin{align*}
&b S' b^{-1}=  S 
, \quad 
(b,e)(S' \times K)(b,e)^{-1} = S \times K \quad \text{and}
\\
&
g P(G_0, S' \times K) g^{-1}= P(G_0, S \times K).
\end{align*}
Thus we can assume $b = e$ without loss of generality. Since $S \subset G_0$ we have $S \cap K = S_{eK}$. Thus by the assumption, for each $a \in U$, $k \in K$ and $h \in P(G_0, U \times K)$ we have 
\begin{align*}
&
a^{-1}S_{aK} a  = S_{eK}
, \quad 
(a,k)^{-1}(S \times K)_{(a,k) \cdot e}(a,k)
=
(S \times K)_{e}
=
\Delta (S_{eK}) \quad \text{and}
\\
&
h^{-1}P(G_0, S \times K)_{h * \hat{0}}h
=
P(G_0, S \times K)_{\hat{0}}
=
\{\hat{s} \in \hat{S} \mid s \in S_{eK}\}.
\end{align*}
Thus by homogeneity it suffices to consider normal vectors only at $eK \in G_0/K$, $e \in G_0$ and $\hat{0} \in V_\mathfrak{g}$. Take $\xi \in T^\perp_{eK}(U \cdot eK)$. Denote by $v \in T^\perp_{e} ((U \times K)\cdot e)$ the horizontal lift of $\xi$ and by $\hat{v} \in T^\perp_{\hat{0}} (P(G_0, U \times K) * \hat{0})$ the horizontal lift of $v$. By similar arguments as in  the proof of Theorem \ref{thm1} we can construct a reflection $\nu_v$ from $\nu_\xi$ and a reflection $\nu_{\hat{\xi}}$ from $\nu_\xi$, and vice versa. These show that the statements (i), (ii) and (iii) in Theorem \ref{thm2} are equivalent and our claim follows.\qed

\section*{Acknoledgements}
The author would like to thank Professors Ernst Heintze and Takashi Sakai for helpful discussions and invaluable suggestions. The author is also grateful to Professor Jost-Hinrich Eschenburg for useful discussions and to Professor Peter Quast for his hospitality during the author's visit to the University of Augsburg. Thanks are also due to Professors  Naoyuki Koike and Hiroshi Tamaru for their interests in his work and valuable comments. Finally the author would like to thank Professor Yoshihiro Ohnita for useful discussions and constant encouragement.

\end{document}